\documentclass[10pt,leqno]{article}
\usepackage{amssymb,amsbsy,amsmath,amsfonts,amssymb,amscd, mathrsfs}
\usepackage{color, graphicx} 

\usepackage[english]{babel}
\usepackage[T1]{fontenc}
\usepackage{indentfirst}

\makeatletter
\@addtoreset{equation}{section}
\makeatother

\newtheorem{statement}{}[section]
\newtheorem{theorem}[statement]{Theorem}
\newtheorem{lemma}[statement]{Lemma} 

\newtheorem{proposition}[statement]{Proposition}

\newcommand\C{\mathbb C}
\newcommand\R{\mathbb R}
\newcommand\T{\mathbb T}
\newcommand\D{\mathbb D}

\newcommand\e{{\rm e}}
\newcommand\eps{\varepsilon}
\newcommand\ind{{\rm 1\kern-.30em I}}
\newcommand\qed{\hfill $\square$}

\let\phi=\varphi

\newcommand\converge{\mathop{\longrightarrow}\limits}

\title{\bf Approximation numbers of composition operators on $H^p$}
\author{\it Daniel Li, Herv\'e Queff\'elec, Luis Rodr{\'\i}guez-Piazza\footnote{Supported by a Spanish research project MTM 2012-05622.}}

\date{\footnotesize \today}

\begin{document}

\maketitle

\noindent{\bf Abstract.} \emph{We give estimates for the approximation numbers of composition operators on the $H^p$ spaces, $1 \leq p < \infty$.}

\medskip

\medskip


\section{Introduction} 

Recently, the study of approximation numbers of composition operators on $H^2$ has been initiated (see \cite{LIQUEROD}, \cite{estimates}, \cite{LLQR}, 
\cite{Queff-Seip}, \cite{spectral radius}), and (upper and lower) estimates have been given. However, most of the techniques used there are specifically 
Hilbertian (in particular Weyl's inequality; see \cite{LIQUEROD}). Here, we consider the case of composition operators on $H^p$ for $1 \leq p < \infty$. 
We focus essentially on lower estimates, because the upper ones are similar, with similar proofs, as in the Hilbertian case. We give in Theorem~\ref{lobo} a 
minoration involving the uniform separation constant of finite sequences in the unit disk and the interpolation constant of their images by the symbol. We finish 
with some upper estimates. 

\subsection{Preliminary}

Recall that if $X$ and $Y$ are two Banach spaces of analytic functions on the unit disk $\D$, and $\phi \colon \D \to \D$ is an analytic self-map of $\D$, one 
says that $\phi$ induces a \emph{composition operator} $C_\phi \colon X \to Y$ if $f \circ \phi \in Y$ for every $f \in X$; $\phi$ is then called the 
\emph{symbol} of the composition operator. One also says that $\phi$ is a symbol for $X$ and $Y$ if it induces a composition operator 
$C_\phi \colon X \to Y$. 
\bigskip

For every $a \in \D$, we denote by $e_a \in (H^p)^\ast$ the evaluation map at $a$, namely:
\begin{equation} 
\qquad e_a (f) = f (a) \, , \quad f \in H^p .
\end{equation} 
We know that (\cite{ZHU}, p.~253):
\begin{equation}\label{eval} 
\| e_a \| = \left(\frac{1}{1 - |a|^2} \right)^{1/p}  
\end{equation}
and the mapping equation 
\begin{equation} 
C_{\varphi}^{\ast} (e_a) = e_{\varphi(a)} 
\end{equation} 
still holds. \par 
\medskip

Throughout this section we denote by $\| \, . \, \|$, without any subscript, the norm in the dual space $(H^p)^\ast$. \par
\smallskip 

Let us stress that this dual norm of $(H^p)^\ast$ is, for $1 < p < \infty$, equivalent, but not equal, to the norm $\| \, . \, \|_q$ of $H^q$, and the 
equivalence constant tends to infinity when $p$ goes to $1$ or to $\infty$. 
\bigskip

As usual, the notation $A \lesssim B$ means that there is a constant $c$ such that $A \leq c \, B$ and $A \approx B$ means that $A \lesssim B$ and 
$B \lesssim A$. 


\subsection{Singular numbers}

For an operator $T \colon X \to Y$ between Banach spaces $X$ and $Y$, its \emph{approximation numbers} are defined, for $n \geq 1$, as:
\begin{equation} \label{approx numbers} 
a_n (T) = \inf_{\text{rank}\, R < n} \| T - R\| \,.
\end{equation} 
One has $\| T \| = a_1 (T) \geq a_2 (T) \geq \cdots \geq a_n (T) \geq a_{n + 1} (T) \geq \cdots$, and (assuming that $Y$ has the Approximation Property), 
$T$ is compact if and only if $a_n (T) \converge_{n \to \infty} 0$.\par\smallskip

We will also need other singular numbers (see \cite{CAST}, p.~49). \par

The \emph{$n$-th Bernstein number} $b_n (T)$ of $T$, defined as:
\begin{equation} \label{Bernstein numbers} 
b_n (T) = \sup_{\substack{E \subseteq X \\ \dim E = n }} \inf_{x \in S_E} \| T x \| \, , 
\end{equation} 
where $S_E = \{ x \in E \, ; \ \| x \| = 1\}$ is the unit sphere of $E$. When these numbers tend to $0$, $T$ is said to be superstrictly singular, or finitely 
strictly singular (see \cite{Plichko}).  \par

The \emph{$n$-th Gelfand number} of $T$, defined as:
\begin{equation} \label{Gelfand numbers} 
c_n (T) = \inf_{\substack{L \subseteq Y \\ {\rm codim}\, L < n}} \, \| T_{\mid L} \| \, ,
\end{equation} 
\smallskip

One always has:
\begin{equation} \label{inegalites} 
\qquad a_n (T) \geq c_n (T) \quad \text{and} \quad a_n (T) \geq b_n (T) \, ,
\end{equation} 
and, when $X$ and $Y$ are Hilbert spaces, one has $a_n (T) = b_n (T) = c_n (T)$ (\cite{Pietsch}, Theorem~2.1). \par


\section{Lower bounds} \label{Lower bounds}

\subsection{Sub-geometrical decay} 

We first show that, as in the Hilbertian case $H^2$ (\cite{LIQUEROD}, Theorem~3.1), the approximation numbers of the composition operators on $H^p$ 
cannot decrease faster than geometrically. 
\par

Though we cannot longer appeal to the Hilbertian techniques of \cite{LIQUEROD}, Weyl's inequality has the following generalization 
(\cite{Carl-Triebel}, Proposition~2). 
\begin{proposition} [Carl-Triebel] \label{Weyl general} 
Let $T$ be a compact operator on a complex Banach space $E$ and $\big(\lambda_n (T) \big)_{n \geq 1}$ be the sequence of its eigenvalues, indexed such that 
$|\lambda_1 (T)| \geq |\lambda_2 (T)| \geq \cdots$. Then, for $n = 1, 2, \ldots$ and $m = 0, 1, \ldots, n - 1$, one has:
\begin{equation} 
\prod_{j = 1}^n |\lambda_j (T)| \leq 16^n \| T \|^m a_{m + 1} (T)^{n - m} \, .
\end{equation} 
\end{proposition}
(see \cite{Carl-Hinrichs} for an optimal result). Then, we can state:
\begin{theorem} 
For every non-constant analytic self-map $\phi \colon \D \to \D$, there exist $0 < r \leq 1$ and $c > 0$, depending only on $\phi$, such that the approximation 
numbers of the composition operator $C_\phi \colon H^p \to H^p$ satisfy:
\begin{displaymath} 
\qquad \quad a_n (C_\phi) \geq c \,  r^n \, , \qquad n = 1, 2, \ldots
\end{displaymath} 
In particular $\liminf_{n \to \infty} [a_n (C_\phi)]^{1/n} \geq r > 0$. 
\end{theorem} 

\noindent{\bf Proof.} If $C_\phi$ is not compact, the result is trivial, with $r = 1$; so we assume that $C_\phi$ is compact. \par\smallskip

Before carrying on, we first recall some notation used in \cite{LIQUEROD}. For every $z \in \D$, let
\begin{displaymath} 
\phi^\sharp (z) = \frac{|\phi ' (z) | \, (1 - |z|^2)}{1 - |\phi (z)|^2} 
\end{displaymath} 
be the pseudo-hyperbolic derivative of $\phi$ at $z$, and 
\begin{displaymath} 
[ \phi ] = \sup_{z \in \D} \phi^\sharp (z) \, .
\end{displaymath} 
By the Schwarz-Pick inequality, one has $[ \phi ] \leq 1$. Moreover, since $\phi$ is not constant, one has $[\phi] > 0$. \par

We also set, for every operator $T \colon H^p \to H^p$:
\begin{displaymath} 
\beta^- (T) = \liminf_{n \to \infty} [a_n (T)]^{1/n} \, .
\end{displaymath} 

For every $a \in \D$, we are going to show that $\beta^- (C_\phi) \geq \big(\phi^\sharp (a) \big)^2$, which will give $\beta^- (C_\phi) \geq [\phi ]^2$, by 
taking the supremum for $a \in \D$, and the stated result, with $0 < r < [\phi]^2$. 
\par\smallskip

If $\phi^\sharp (a) = 0$, the result is obvious, so we assume that $\phi^\sharp (a) > 0$. \par

We consider the automorphism $\Phi_a$, defined by $\Phi_a (z) = \frac{a - z}{1 - \overline{a} z}$, and set
\begin{displaymath} 
\psi_a = \Phi_{\phi (a)} \circ \phi \circ \Phi_a \, . 
\end{displaymath} 
One has $\psi_a (0) = 0$ and $|\psi_a ' (0) | = \phi^\sharp (a)$. \par\smallskip

Since $C_\phi$ is compact on $H^p$, $C_{\psi_a} = C_{\Phi_a} \circ C_\phi \circ C_{\Phi_{\phi (a)}}$ is also compact on $H^p$. But we know that this 
is equivalent to say that it is compact on $H^2$. Since $\psi_a (0) = 0$ and $\psi_a ' (0) = \phi^\sharp (a) \neq 0$, we know, by the Eigenfunction Theorem 
(\cite{Shapiro-livre}, p.~94), that the eigenvalues of $C_{\psi_a} \colon H^2 \to H^2$ are the numbers $\big( \psi_a ' (0) \big)^j$, $j = 0, 1, \ldots$, and have 
multiplicity one.  Moreover, the proof given in \cite{Shapiro-livre}, \S~6.2 shows that the eigenfunctions $\sigma^j$ are not only in $H^2$, but in all $H^q$, 
$1 \leq q < \infty$. Hence $\lambda_j (C_{\psi_a}) = \big( \psi_a ' (0) \big)^{j - 1}$. We now use Proposition~\ref{Weyl general}, with $2n$ instead of $n$ and 
$m = n - 1$; we get:
\begin{align*} 
|\psi_a ' (0) |^{n (2 n - 1)} 
& = \prod_{j = 1}^{2 n} |\lambda_j (C_{\psi_a}) | \leq 16^{2 n} \| C_{\psi_a} \|^{n - 1}  a_n (C_{\psi_a})^{n + 1}  \\ 
& \leq 16^{2 n} \| C_{\psi_a} \|^n  a_n (C_{\psi_a})^n \, , 
\end{align*} 
since $a_n (C_{\psi_a} ) \leq \| C_{\psi_a} \|$.  \par

That implies that $\beta^- (C_{\psi_a}) \geq |\psi_a ' (0)|^2 = \big( \phi^\sharp (a) \big)^2$. \par

Since $C_{\Phi_a}$ and $C_{\Phi_{\phi (a)}}$ are automorphisms, we have $\beta^- (C_\phi) = \beta^- (C_{\psi_a}) $, hence the result. \qed


\subsection{Main result} 

In this section, we use the fortunate fact that, though the evaluation maps at well-chosen points of $\D$ can no longer be said to constitute a Riesz sequence, they 
will still constitute an unconditional sequence in $H^p$ with good constants, as we are going to see, which will be sufficient for our purposes. \par
\medskip

Recall (see \cite{GAR}, p.~276) that the \emph{interpolation constant} $\kappa_\sigma$ of a finite sequence $\sigma = (z_1, \ldots, z_n)$ of points 
$z_1, \ldots, z_n \in \D$ is defined by:
\begin{equation} 
\kappa_\sigma = \sup_{|a_1|, \, \ldots, |a_n| \leq 1} \inf\{ \| f \|_\infty \, ; \ f \in H^\infty \text{ and } f (z_j) = a_j \, , 1 \leq j \leq n \} \, .
\end{equation} 
Then:
\begin{lemma} \label{same} 
For \hfil every \hfil finite \hfil sequence \hfil $\sigma = (z_1, \ldots, z_n)$ \hfil of \hfil distinct \hfil points \hfil $z_1, \ldots, z_n \in \D$, one has:
\begin{equation} \label{observ} 
\kappa_\sigma^{- 1}\Big\| \sum_{j = 1}^n \lambda_j e_{z_j} \Big\|  
\leq \Big\| \sum_{j = 1}^n \omega_j \lambda_j e_{z_j}\Big \| 
\leq \kappa_\sigma \Big\| \sum_{j = 1}^n \lambda_j e_{z_j} \Big\| 
\end{equation} 
for all $\lambda_1, \ldots, \lambda_n \in \C$ and all complex numbers numbers $\omega_1, \ldots, \omega_n$ such that 
$|\omega_1| =  \cdots = |\omega_n| = 1$.
\end{lemma}

\noindent{\bf Proof.} Set $L = \sum_{j = 1}^n \lambda_j e_{z_j}$ and $L_\omega = \sum_{j = 1}^n \omega_j \lambda_j e_{z_j}$. There exists 
$h \in H^\infty$ such that $\| h \|_\infty \leq \kappa_\sigma$ and $h ( z_j) = \omega_j$ for every $j = 1, \ldots, n$. For every $g \in H^p$, one has 
$L_\omega (g) = \sum_{j = 1}^n \omega_j \lambda_j g (z_j) = \sum_{j = 1}^n h (z_j) \lambda_j g (z_j) = L (h g)$; hence:
\begin{displaymath} 
| L_\omega (g) | \leq \| L \| \, \| h g \|_p \leq \| L \| \, \| h \|_\infty \| g \|_p \leq \kappa_\sigma \| L \| \, \| g \|_p
\end{displaymath} 
and we get $\| L_\omega \| \leq \kappa_\sigma \| L \|$, which is the right-hand side of \eqref{observ}. The left-hand side follows, by replacing 
$\lambda_1, \ldots, \lambda_n$ by $\overline{\omega_1} \lambda_1, \ldots, \overline{\omega_n} \lambda_n$. \qed
\par\medskip

We now prove the following lower estimate. 
\begin{theorem} \label{lobo} 
Let $\phi \colon \D \to \D$ and $C_\phi \colon H^p \to H^p$, with $1 \leq p < \infty$. Let $u_1, \ldots, u_n \in \D$ such that 
$v_1 =\phi (u_1), \ldots, v_n = \phi (u_n)$ are distinct. Then, for some constant $c_p$ depending only on $p$, we have:
\begin{equation} \label{zero} 
a_{n} (C_\phi) \geq c_p \, 
\kappa_v^{- 1} \Big(1 + \log \frac{1}{\delta_u} \Big)^{- 1/\min (p, 2)} \inf_{1 \leq j \leq n} \bigg(\frac{1 - |u_j|^2}{1 - |v_j|^2} \bigg)^{1/p} 
\, ,
\end{equation} 
where $\delta_u$ is the uniform separation constant of the sequence $u = (u_1, \ldots, u_n)$  and $\kappa_v$ the interpolation constant of 
$v = (v_1, \ldots, v_n)$.  
\end{theorem} 

For the proof, we need to know some precisions on the constant in Carleson's embedding theorem. Recall that the \emph{uniform separation constant} 
$\delta_\sigma$ of a finite sequence $\sigma = (z_1, \ldots, z_n)$ in the unit disk $\D$, is defined by:
\begin{equation} \label{Carleson constant}
\delta_\sigma = \inf_{1 \leq j \leq n} \prod_{k \neq j} \Big| \frac{z_j - z_k}{1 - \overline{z_j} z_k} \Big| \, \cdot
\end{equation}
\begin{lemma}  \label{Carleson} 
Let $\sigma = (z_1, \ldots, z_n)$ be a finite sequence of distinct points in $\D$ with uniform separation constant $\delta_\sigma$. Then:
\begin{equation} 
\sum_{j = 1}^n (1 - |z_j|^2) \, |f (z_j)|^p  \leq 12 \, \bigg[ 1 + \log \frac{1}{\delta_\sigma} \bigg] \, \| f \|_p^p 
\end{equation} 
for all $f \in H^p$. 
\end{lemma} 

\noindent{\bf Proof.} For $a \in \D$, let $k_a (z) = \frac{\sqrt{1 - |a|^2}}{1 - \overline{a} z}$ be the normalized reproducing kernel. For every positive 
Borel measure $\mu$ on $\D$, let:
\begin{displaymath} 
\gamma_\mu = \sup_{a \in {\rm supp}\, \mu} \int_\D | k_a (z)|^2 \, d\mu (z) \, .
\end{displaymath} 
The so-called Reproducing Kernel Thesis (see \cite{NIK}, Lecture~VII, pp.~151--158) says that there is an absolute positive constant $A_1$ such that:
\begin{displaymath} 
\int_\D | f (z) |^p \, d \mu (z) \leq A_1 \, \gamma_\mu \, \| f  \|_p^p 
\end{displaymath} 
for every $f \in H^p$ (that follows from the case $p = 2$ in writing $f = B h^{2/p}$ where $B$ is a Blaschke product and $h \in H^2$). Actually, one can take 
$A_1 = 2 \, \e$ (see \cite{Steph}, Theorem~0.2). But when $\mu$ is the discrete measure $\sum_{j = 1}^n (1 - |z_j|^2) \, \delta_{z_j}$, it is not difficult to 
check (see \cite{Duren}, Lemma~1, p.~150, or \cite{HOF}, p.~201) that:
\begin{displaymath} 
\gamma_\mu \leq 1 + 2 \, \log \frac{1}{\delta_\sigma} \, \cdot
\end{displaymath} 
That gives the result since $4 \, \e \leq 12$. \qed 
\par\medskip 

\noindent{\bf Proof of Theorem~\ref{lobo}.} We will actually work with the Bernstein numbers of $C_\phi^\ast$. Recall that they are defined
in \eqref{Bernstein numbers}. That will suffice since $a_n (C_\phi) \geq a_n (C_\phi^\ast )$ (one has equality if $C_\phi$ is compact: see \cite{Victoria} or 
\cite{CAST}, pp.~89--91) and $a_n (C_\phi^\ast )\geq b_n (C_\phi^\ast)$. \par
\smallskip 

Take $u_1, \ldots, u_n \in \D$ such that $v_1 = \phi (u_1), \ldots, v_n = \phi (u_n)$ are distinct. The points $u_1, \ldots, u_n$ are then also distinct and the 
subspace $E = {\rm span} \, \{e_{u_1}, \ldots, e_{u_n} \}$ of $(H^p)^\ast$ is $n$-dimensional. Let 
\begin{displaymath} 
L = \sum_{j = 1}^n \lambda_j e_{u_j} 
\end{displaymath} 
be in the unit sphere of $E$. We set, for $f \in H^p$ and for $ j = 1, \ldots, n$:
\begin{displaymath} 
\qquad \Lambda_j = \lambda_j \, \| e_{u_j} \| \, , \quad \text{and} \quad F_j = \| e_{u_j} \|^{- 1} f (u_j) \, , 
\end{displaymath} 
and finally: 
\begin{displaymath} 
\Lambda = (\Lambda_1, \ldots, \Lambda_n) \quad \text{and} \quad  F = (F_1, \ldots, F_n) \, . 
\end{displaymath} 
\par\smallskip

We will separate three cases. \par
\smallskip

\noindent {\bf Case 1:} $1 < p \leq 2$. \par

One has $\| C_{\phi}^{\ast} (L) \|  = \big\| \sum_{j = 1}^n \lambda_j \, e_{v_j} \big\|$. Using Lemma~\ref{same}, we 
obtain for any choice of complex signs $\omega_1, \ldots,  \omega_n$:
\begin{equation} \label{a moyenner} 
\| C_{\phi}^{\ast} (L) \| \geq \kappa_v^{- 1}\Big\| \sum_{j = 1}^n \omega_j \lambda_j e_{v_j} \Big\| \, . 
\end{equation} 
Let now $q$ be the conjugate exponent of $p$. We know that the space $H^p$ is of type $p$ as a subspace of $L^p$ (\cite{LQ}, p.~169) and therefore its dual 
$(H^p)^\ast$ is of cotype $q$ (\cite{LQ}, p.~165), with cotype constant $\leq \tau_p$, the type $p$ constant of $L^p$ (let us note that we might use that 
$(H^p)^\ast$ is isomorphic to the subspace $H^q$ of $L^q$, but we have then to introduce the constant of this isomorphism). Hence, by averaging 
\eqref{a moyenner} over all independent choices of signs and using the cotype $q$ property of $(H^p)^\ast$, we get:
\begin{displaymath} 
\| C_\phi^\ast (L) \|  
\geq \tau_p^{- 1} \, \kappa_v^{- 1} \Big( \sum_{j = 1}^n | \lambda_j |^q \| e_{v_j} \|^q \Big)^{1/q} \notag \\ 
\geq \tau_p^{- 1} \, \kappa_v^{- 1} \mu_n \Big( \sum_{j = 1}^n |\lambda_j|^q \| e_{u_j} \|^q \Big)^{1/q} \, ,
\end{displaymath}
so that  
\begin{equation} \label{etoile} 
\| C_\phi^\ast (L) \|  \geq \tau_p^{- 1} \, \kappa_v^{- 1} \mu_n \, \| \Lambda \|_q \, , 
\end{equation} 
where:
\begin{displaymath} 
\mu_n = \inf_{1 \leq j \leq n} \frac{\| e_{v_j} \|}{\| e_{u_j} \|} = \inf_{1 \leq j \leq n} \Big( \frac{1 - |u_j|^2}{1 - |v_j|^2} \Big)^{1/p} \, .
\end{displaymath} 
\par\smallskip

It remains to give a lower bound for $\| \Lambda \|_q$. \par

But, by H\"older's inequality:
\begin{displaymath} 
| L (f) | = \Big| \sum_{j = 1}^n \lambda_j f (u_j) \Big| = \Big| \sum_{j = 1}^n \Lambda_j F_j \Big| 
\leq \| \Lambda \|_q \| F \|_p \, . 
\end{displaymath} 
Since 
\begin{displaymath} 
\| F \|_p^p = \sum_{j = 1}^n \| e_{u_j} \|^{- p} | f (u_j)|^p = \sum_{j = 1}^n (1 - |u_j|^2) \, | f (u_j)|^p , 
\end{displaymath} 
Lemma~\ref{Carleson} gives:
\begin{displaymath} 
| L (f) | \leq \| \Lambda \|_q \bigg[ 12 \, \Big( 1 + \log \frac{1}{\delta_u} \Big) \bigg]^{1/p} \| f \|_p .
\end{displaymath} 
Taking the supremum over all $f$ with $\| f \|_p \leq 1$, we get, taking into account that $\| L \| = 1$:
\begin{equation} \label{lambda} 
\| \Lambda \|_q \geq  \bigg[ 12 \, \Big( 1 + \log \frac{1}{\delta_u} \Big) \bigg]^{- 1/p} \, . 
\end{equation} 

By combining \eqref{etoile} and \eqref{lambda}, we get:
\begin{displaymath} 
\| C_\phi^\ast (L) \| \geq  (12)^{- 1/p} \, \tau_p^{- 1} \, \mu_n \, \kappa_v^{ - 1} \Big( 1 + \log \frac{1}{\delta_u} \Big)^{- 1/p} \, .
\end{displaymath} 

Therefore:
\begin{displaymath} 
b_n (C_\phi^\ast) \geq (12)^{- 1/p} \, \tau_p^{- 1} \, \mu_n \, \kappa_v^{ - 1} \Big( 1 + \log \frac{1}{\delta_u} \Big)^{- 1/p} \, .
\end{displaymath} 
\par\smallskip 

\noindent{\bf Case 2:} $2 < p < \infty$. \par

We follow the same route, but in this case, $H^p$ is of type $2$ and hence $(H^p)^\ast$ is of cotype $2$. Therefore, we get:
\begin{equation} 
\| C_\phi^\ast (L) \| \geq \tau_2^{- 1} \, \kappa_v^{- 1} \, \mu_n \, \| \Lambda \|_2 
\end{equation} 
and, using Cauchy-Schwarz inequality: 
\begin{equation} 
\| \Lambda \|_2 \geq \bigg[ 12 \Big( 1 + \log \frac{1}{\delta_u} \Big) \bigg]^{- 1/2} \, ; 
\end{equation} 
so:
\begin{equation} 
\| C_\phi^\ast (L) \| \geq (12)^{- 1/2} \, \tau_2^{- 1} \, \mu_n \, \kappa_v^{- 1} \Big( 1 + \log \frac{1}{\delta_u} \Big)^{- 1/2} \, .
\end{equation} 
\par\smallskip 

\noindent{\bf Case 3:} $p = 1$. \par

In this case $(H^1)^\ast$ (which is isomorphic to the space $BMOA$) has no finite cotype. But, for each $k = 1, \ldots, n$, one has, using Lemma~\ref{same}:
\begin{align*} 
|\lambda_k| \, \| e_{v_k} \| 
& = \frac{1}{2} \, \bigg\| \bigg(\sum_{j \neq k} \lambda_j e_{v_j} + \lambda_k e_{v_k} \bigg) 
- \bigg(\sum_{j \neq k} \lambda_j e_{v_j} - \lambda_k e_{v_k} \bigg) \bigg\| \\ 
& \leq \frac{1}{2} \, \bigg( \bigg\| \sum_{j \neq k} \lambda_j e_{v_j} + \lambda_k e_{v_k} \bigg\| 
+ \bigg\| \sum_{j \neq k} \lambda_j e_{v_j} - \lambda_k e_{v_k} \bigg\|  \bigg) \\ 
& \leq \kappa_v \bigg\| \sum_{j = 1}^n \lambda_j e_{v_j} \bigg\| \, ;
\end{align*} 
hence:
\begin{equation} 
\| C_\phi^\ast (L) \| \geq \kappa_v^{- 1} \, \mu_n \, \| \Lambda \|_\infty \, . 
\end{equation} 
Since $| L (F) | \leq \| \Lambda \|_\infty \| F \|_1$, we get, as above, using Lemma~\ref{Carleson}:
\begin{equation} 
\| \Lambda \|_\infty \geq \bigg[ 12 \, \Big( 1 + \log \frac{1}{\delta_u} \Big) \bigg]^{- 1} \, , 
\end{equation} 
and therefore:
\begin{equation} 
\| C_\phi^\ast (L) \| \geq (12)^{- 1} \, \mu_n \, \kappa_v^{ - 1} \, \Big( 1 + \log \frac{1}{\delta_u} \Big)^{- 1} 
\end{equation} 
and that finishes the proof of Theorem~\ref{lobo}. \qed 
\par
\bigskip\goodbreak

\noindent {\bf Example.} We will now apply this result to lens maps. We refer to \cite{Shapiro-livre} or \cite{LLQR} for their definition. For $\theta \in (0, 1)$, 
we denote:
\begin{equation} 
\lambda_\theta (z) = \frac{(1 + z)^\theta - (1 - z)^\theta}{(1 + z)^\theta + (1 - z)^\theta} \, \cdot
\end{equation} 
\begin{proposition} \label{appl} 
Let $\lambda_\theta$ be the lens map of parameter $\theta$ acting on $H^p$, with $1 \leq p < \infty$. Then, for positive constants $a$ and $b$, depending 
only on $\theta$ and $p$: 
\begin{displaymath} 
a_n (C_{\lambda_\theta}) \geq a \, \e^{- b \sqrt n} . 
\end{displaymath} 
\end{proposition} 

Actually, this estimate is valid for polygonal maps as well. 
\par\medskip

\noindent {\bf Proof.} Let $0 < \sigma < 1$ and consider $u_j = 1 - \sigma^j$ and $v_j = \lambda_\theta (u_j)$, $1 \leq j \leq n$. We know 
from \cite{LIQUEROD}, Lemma~6.4 and Lemma~6.5, that, for $\alpha = \frac{\pi^2}{2}$ and $\beta = \beta_\theta = \frac{\pi^2}{2^\theta \theta}$:
\begin{displaymath} 
\delta_u  \geq \e^{- \alpha / (1 - \sigma)}  \quad \text{and} \quad \delta_v \geq \e^{- \beta / (1 - \sigma)} .
\end{displaymath} 

But we know that the interpolation constant $\kappa_\sigma$ is related to the uniform separation constant $\delta_\sigma$ by the following inequality 
(\cite{GAR} page~278), in which $\Lambda$ is a positive numerical constant:
\begin{equation} \label{Interpolation constant}
\frac{1}{\delta_\sigma} \leq \kappa_\sigma \leq \frac{\Lambda}{\delta_\sigma} \bigg(1 + \log \frac{1}{\delta_\sigma}\bigg) \cdot
\end{equation}
Actually, S. A. Vinogradov, E. A. Gorin and S. V. Hru{\v s}c\"ev \cite{Vino} (see \cite{Mortini}, p.~505) proved that 
\begin{displaymath} 
\kappa_\sigma \leq \frac{2 \, \e}{\delta_\sigma} \, \Big( 1 + 2 \log \frac{1}{\delta_\sigma} \Big) \, \raise 1pt \hbox{,} 
\end{displaymath} 
so we can take $\Lambda \leq 4 \, \e \leq 12$. \par

It follows that 
\begin{equation} \label{morceau 1}
\kappa_v^{- 1} \geq \frac{1 - \sigma}{\Lambda (\beta + 1)} \, \e^{- \beta / (1 - \sigma)}. 
\end{equation} 

Setting $\tilde p = \min (p, 2)$, we have:
\begin{equation} \label{morceau 2}
\Big(1 + \log \frac{1}{\delta_u} \Big)^{- 1 / \tilde p} \geq \Big( \frac{1 - \sigma}{\alpha + 1} \Big)^{1/ \tilde p} \, .
\end{equation} 

We now estimate $\mu_n$. \par
Since $\lambda_\theta (0) = 0$, Schwarz's lemma says that $|\lambda_\theta (z) | \leq |z|$; hence 
$\frac{1 - |z|^2}{1 - |\lambda_\theta (z)|^2} \geq \frac{1 - |z|}{1 - |\lambda_\theta (z)|}$. 
But $1 - v_j = 1 - \lambda_\theta (u_j) = \frac{2 \sigma^{j \theta}}{(2 - \sigma^j)^\theta + \sigma^{j \theta}}$; hence (since $u_j$ and $v_j$ are real):
\begin{displaymath} 
\frac{1 - |u_j|^2}{1 - |v_j|^2} \geq \frac{1 - u_j}{1 - v_j} = \frac{\sigma^j}{2 \sigma^{j \theta}} \, [(2 - \sigma^j)^\theta + \sigma^{j \theta}] \,.
\end{displaymath} 
Since the function $f (x) = (2 - x)^\theta + x^\theta$ increases on $[0, 1]$, one gets:
\begin{displaymath} 
\frac{1 - |u_j|^2}{1 - |v_j|^2} \geq \Big( \frac{1}{2} \, \sigma^j \Big)^{1 - \theta} \, ,
\end{displaymath} 
and therefore:
\begin{equation} \label{morceau 3}
\mu_n \geq \Big( \frac{1}{2} \, \sigma^n \Big)^{(1 - \theta) / p} \, .
\end{equation} 

Applying now Theorem~\ref{lobo} and using \eqref{morceau 1}, \eqref{morceau 2} and \eqref{morceau 3}, we get: 
\begin{displaymath} 
a_n (C_{\lambda_\theta}) \geq  \alpha_{p, \theta} \, \e^{- \beta / (1 - \sigma)}\, (1 - \sigma)^{1/ \tilde p} \, \sigma^{n (1 - \theta) / p} 
\end{displaymath} 
with $\alpha_{p, \theta} = \frac{c_p}{\Lambda (\beta + 1) (\alpha + 1)^{1/ \tilde p} 2^{(1 - \theta)/p}}\, \cdot$ \par

Taking $\sigma = \e^{- \varepsilon}$ where $0 < \varepsilon < 1$, we get, since $1 - \e^{ - \eps} \geq \eps /2$: 
\begin{displaymath} 
a_n (C_{\lambda_\theta}) \geq  \alpha_{p, \theta} \, \e^{- 2 \beta / \eps} \, \Big(\frac{\eps}{2}\Big)^{1/ \tilde p} \,  \e^{- \eps n (1 - \theta)/ p}. 
\end{displaymath} 
Optimizing by taking $\varepsilon = \sqrt{\frac{3 \beta p}{1 - \theta}} \, \frac{1}{\sqrt n}$ gives, for $n$ large enough (in order to have $\eps < 1$): 
\begin{equation} \label{lens}
a_n (C_{\lambda_\theta}) \geq \alpha'_{p, \theta} \, n^{- 1 / (2 \tilde p)} \, \e^{ - \beta_{p, \theta} \sqrt{n}} 
\end{equation} 
with $\alpha'_{p, \theta} = \alpha_{p, \theta} \big( \frac{\beta p}{2 (1 - \theta)} \big)^{1/ (2 \tilde p)}$ and 
$\beta_{p, \theta} = \sqrt{ \frac{2 \beta (1 - \theta)}{p}}\, \cdot$ \par\smallskip

We get Theorem~\ref{appl}, with $b > \beta_{p, \theta}$. \qed
\par\medskip

Let us note that $\beta_{p, \theta} = \frac{2^{\frac{1 - \theta}{2}} \pi}{\sqrt{p}} \, \sqrt{\frac{1 - \theta}{\theta}}$ tends to $0$ when $\theta$ goes to $1$ 
and tends to infinity when $\theta$ goes to $0$. 
\goodbreak


\subsection{A minoration depending on the radial behaviour of $\phi$} 

We are using Theorem~\ref{lobo} to give, as in \cite{estimates}, Theorem~3.2, a lower bound for $a_n (C_\phi)$ which depends on the behaviour of $\phi$ 
near $\partial \D$. \par
\medskip

We recall first (see \cite{estimates}, Section~3) that an analytic self-map $\phi \colon \D \to \D$ is said to be \emph{real} if it takes real values on $]-1, 1[$. 
If $\omega \colon [0, 1] \to [0, 2]$ is a modulus of continuity (meaning that $\omega$ is continuous, increasing, sub-additive, vanishing at $0$, and concave), 
$\phi$ is said to be an \emph{$\omega$-radial symbol} if it is real and:
\begin{equation} \label{radial symbol}
\qquad 1 - \phi (r) \leq \omega (1 - r)\, , \quad 0 \leq r < 1 \, .
\end{equation} 
We have the following result. \goodbreak

\begin{theorem} \label{mino radial} 
Let $\phi$ be an $\omega$-radial symbol. Then, for $1 \leq p < \infty$, the approximation numbers of the composition operator $C_\phi \colon H^p \to H^p$ 
satisfy:
\begin{equation} \label{eq mino radial}
a_n (C_\phi) \geq  
c_p' \sup_{0 < \sigma < 1} \! \bigg[ \Big(\frac{\omega^{- 1} (a \, \sigma^n)}{a \, \sigma^n} \Big)^{1/p} 
( 1 - \sigma)^{1 / \max (p^\ast, 2)} \,  \exp \Big( \! - \frac{5}{1 - \sigma} \Big) \bigg] , 
\end{equation} 
where $c'_p$ is a constant depending only on $p$, $p^\ast$ is the conjugate exponent of $p$, and $a = 1 - \phi (0) > 0$. 
\end{theorem} 

\noindent{\bf Proof.} As in \cite{estimates}, p.~556, we fix $0 < \sigma < 1$ and define inductively $u_j \in [0, 1)$ by $u_0 = 0$ and, using the intermediate 
value theorem:
\begin{displaymath} 
\qquad 1 - \phi (u_{j + 1}) = \sigma \, [1 - \phi (u_j)] \, , \quad \text{with } 1 > u_{j + 1} > u_j \, .
\end{displaymath} 
We set $v_j = \phi (u_j)$. We have $- 1 < v_j < 1$ and $1 - v_n = a\, \sigma^n$. We proved in \cite{estimates}, p.~556, that:
\begin{equation} \label{one} 
\frac{1 - |u_j|^2}{1 - |v_j|^2} \geq \frac{1}{2}\, \frac{ \omega^{- 1} (a\, \sigma^n)}{a \, \sigma^n} \, \cdot
\end{equation} 
Moreover, we proved in \cite{estimates}, p.~557, that the uniform separation constant of $v = (v_1, \ldots, v_n)$ is such that:
\begin{equation} \label{two}
\delta_v \geq \exp \Big( - \frac{5}{1 - \sigma} \Big) \, \cdot
\end{equation} 
Since $\delta_u \geq \delta_v$, we get, from \eqref{Interpolation constant}, that:
\begin{equation} \label{three}
\kappa_u \leq 12 \, \Big( \frac{6 - \sigma}{1 - \sigma} \Big)\, \exp\Big( \frac{5}{1 - \sigma} \Big)  
\leq 60 \, \, \Big( \frac{1}{1 - \sigma} \Big)\, \exp\Big( \frac{5}{1 - \sigma} \Big)  \, \cdot
\end{equation} 
Using now \eqref{zero} of Theorem~\ref{lobo} and combining \eqref{one}, \eqref{two} and \eqref{three}, we get Theorem~\ref{mino radial}. \qed
\par\bigskip\goodbreak

\noindent{\bf Example 1: lens maps.} Let us come back to the lens maps $\lambda_\theta$ for testing Theorem~\ref{mino radial}. We have 
$\omega^{- 1} (h) \approx h^{1/\theta}$ (see \cite{LLQR}, Lemma~2.5) and $ a = 1 - \lambda_\theta (0) = 1$. Setting 
$K = \frac{1}{10 \sqrt{p}} \, \sqrt{\frac{1 - \theta}{\theta}}$ and taking, for $n$ large enough, $\sigma = 1 - \frac{1}{K \sqrt{n}}$, we have, using that 
$\e^{- s} \leq 1 - \frac{4}{5} s $ for $s > 0$ small enough, $\sigma^n \geq \exp ( - \frac{5}{4 K}\, \sqrt{n})$ and hence: 
\begin{displaymath} 
a_n (C_{\lambda_\theta}) \geq c_{\theta, p} \, n^{- \frac{1}{2 \, \max (p^\ast, 2)}} \, 
\exp \bigg[ - \frac{5}{\sqrt{p}} \, \sqrt{\frac{1 - \theta}{\theta}} \, \sqrt{n} \bigg] \, .
\end{displaymath} 
\par

Note that the coefficient of $\sqrt n$ in the exponential is slightly different of that in \eqref{lens}, but of the same order. \par
\medskip

\begin{figure}[ht]
\centering
\includegraphics[width=5cm]{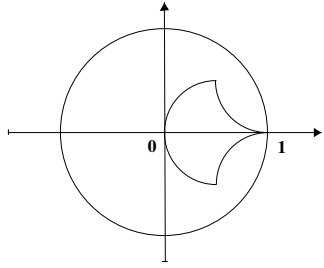}
\caption{\it Cusp map domain} \label{cusp map}
\end{figure}

\noindent{\bf Example 2: cusp map.} We refer to \cite{estimates}, Section~4, for its definition and properties. It is the conformal mapping $\chi$ from $\D$ 
onto the domain represented on Fig.~\ref{cusp map} such that  $\chi (1) = 1$, $\chi (- 1) = 0$, $\chi (i) = (1 + i)/2$ and $\chi (- i) = (1 - i)/2$.
We proved in \cite{estimates}, Lemma~4.2, that, for $0 \leq r < 1$, one has:
\begin{displaymath} 
1 - \chi (r) = \frac{1}{1 + \frac{2}{\pi} \, \log \big[ 1 / 2\,  \arctan \big( \frac{1 - r}{1 + r} \big) \big] } \, \cdot
\end{displaymath} 
Since $1 - \frac{2}{\pi} \log 2 > 0$ and $\arctan x \leq x$ for $x \geq 0$, we get that:
\begin{displaymath} 
1 - \chi (r) \leq \frac{\pi}{2}\, \frac{1}{\log \big( \frac{1 + r}{1 - r} \big)} \leq \frac{\pi}{2}\, \frac{1}{\log \big( \frac{1}{1 - r} \big)} 
\leq 2 \, \, \frac{1}{\log \big( \frac{1}{1 - r} \big)} \cdot
\end{displaymath} 
Hence $\chi$ is an $\omega$-radial symbol with $\omega (x) = 2 / \log (1/x)$. Then $\omega^{ - 1} (h) = \e^{- 2 /h}$. By choosing 
$\sigma = 1 - \frac{\log n}{4 n}$ in \eqref{eq mino radial}, we get, using that $\log (1 - x) \geq - 2 x$ for $x > 0$ small enough, that, for $n$ large enough, 
$\sigma^n \geq 1 / \sqrt{n}$; hence:
\begin{displaymath} 
a_n (C_\chi) \geq c''_p \Big( \sqrt{n} \, \exp \big[ - (2\ a) \, \sqrt{n} \big] \Big)^{1/p} \, \Big(\frac{\log n}{n} \Big)^{1/ \max (p^\ast, 2)} 
\, \exp \Big( - \frac{20 n}{\log n} \Big) \, \cdot
\end{displaymath} 
It follows that, for some constant $C_p > 0$ depending only on $p$, we have:
\begin{equation} 
a_n (C_\chi) \geq C_p \, \exp \Big( - \frac{ 25 n}{\log n} \Big) \, \cdot
\end{equation} 
It has to be stressed that the term in the exponential does not depend on $p$. 
\par\medskip\goodbreak

\noindent{\bf Example 3: Shapiro-Taylor's maps.} 
These maps $\varsigma_\theta$, for $\theta > 0$, were defined in \cite{Shap-Taylor}. Let us recall their definition. For $\eps > 0$, we set 
$V_\eps = \{ z\in \C \, ; \ \Re z > 0 \text{ and } |z | < \eps \}$. For $\eps = \eps_\theta > 0$ small enough, one can define 
\begin{equation} 
f_\theta (z) = z (- \log z )^\theta ,
\end{equation} 
for $z \in V_\eps$, where $\log z$ will be the principal determination of the logarithm. Let now $g_\theta$ be the conformal mapping from $\D$ onto $V_\eps$, 
which maps $\T = \partial \D$ onto $\partial V_\eps$, defined by $g_\theta (z) = \eps\, \phi_0 (z)$, where $\phi_0$ is the conformal map from $\D$ onto 
$V_1$, given by:
\begin{equation} 
\phi_0 (z) = \frac{\displaystyle \Big( \frac{z - i}{i z - 1} \Big)^{1/2} - i} {\displaystyle - i \, \Big( \frac{z - i}{i z - 1} \Big)^{1/2} + 1} \, \cdot
\end{equation} 
Then, we define:
\begin{equation} 
\varsigma_\theta = \exp ( - f_\theta \circ g_\theta) .
\end{equation} 

We saw in \cite{estimates}, p.~560, that 
$\omega^{ - 1} (h) = K_\theta \, h \big( \log (1/h) \big)^{- \theta}$. Hence, choosing $\sigma = 1 / (\e \, \alpha_\theta^{1/n})$, where 
$\alpha_\theta = 1 - \varsigma_\theta (0)$, we get that:
\begin{equation} 
a_n (C_{\varsigma_\theta}) \geq c_{p, \theta} . \frac{1}{n^{\theta/ 2 p}} \, \cdot
\end{equation} 
However, we already remarked in \cite{estimates}, Section~4.2, that, even for $p = 2$, this result is not optimal. 
\goodbreak
 
\section{Upper bound} 

For upper bounds, there is essentially no change with regard to the case $p = 2$. Hence we essentially only state some results. \par\smallskip

We have the following upper bound, which can be obtained with the same proof as in \cite{LLQR}. 
\begin{theorem} \label{upbo} 
Let  $C_\phi \colon H^p \to H^p$, $1 \leq p < \infty$, a composition operator, and $n \geq 1$. Then, for every Blaschke product $B$ with (strictly) less than $n$ 
zeros, each counted with its multiplicity, one has:
\begin{displaymath} 
a_n (C_\phi)\leq C \sqrt {n} \, \Bigg( \sup_{\substack{0 < h < 1 \\ \xi \in \T}} \frac{1}{h} \int_{\overline{S (\xi, h)}} |B|^p \, dm_{\phi} \Bigg)^{1/p} ,
\end{displaymath} 
where $m_\varphi$ is the pullback measure of $m$, the normalized Lebesgue measure on $\T$, under $\phi$ and $S (\xi, h) = \D \cap D (\xi, h)$ is the 
Carleson window of size $h$ centered at $\xi \in \T$. 
\end{theorem} 

\noindent{\bf Proof.} We first estimate the Gelfand number $c_n (C_\phi)$ by restricting to the subspace $B H^p$ which is of codimension $ < n$. 
As in \cite{LLQR}, Lemma~2.4:
\begin{displaymath} 
c_n (C_\phi) \lesssim \Bigg( \sup_{\substack{0 < h < 1 \\ \xi \in \T}} \frac{1}{h} \int_{ \overline{S (\xi, h)}} |B|^p \, dm_\phi \Bigg)^{1/p} .
\end{displaymath} 

Now (see \cite{CAST}, Proposition~2.4.3), one has $a_n (C_\phi) \leq \sqrt{2 n} \, c_n (C_\phi)$, hence the result. \qed
\par\medskip 

We can then deduce, with the same proof, the following version of \cite{estimates}, Theorem~2.3. \par\smallskip

Recall (\cite{estimates}, Definition~2.2) that a symbol $\phi \in A (\D)$ (i.e. $\phi \colon \overline{\D} \to \overline{\D}$ is continuous and analytic in $\D$) 
is said to be \emph{globally regular} if  $\phi (\overline \D) \cap \partial \D = \{\xi_1, \ldots, \xi_l \}$ and there exists a modulus of continuity $\omega$ 
(i.e. a continuous, increasing and sub-additive function $\omega \colon [0, A] \to \R^+$, which vanishes at zero, and that we may assume to be concave), such 
that, writing $E_{\xi_j} = \{ t \, ; \ \gamma (t) = \xi_j \}$, one has $\T = \bigcup_{j = 1}^l \big( E_{\xi_j} + [- r_j, r_j] \big)$ for some $r_1, \ldots, r_l > 0$, 
and for some positive constants $C, c > 0$: 
\begin{align} 
|\gamma (t) - \gamma (t_j) | & \leq C \big(1 - | \gamma (t) | \big) \label{reluc glob} \\
c \, \, \omega (| t - t_j |) & \leq | \gamma (t) - \gamma (t_j) |  \label{ponct glob}
\end{align}
for $j = 1, \ldots, l$, all $t_j \in E_{\xi_j}$ with $|t - t_j| \leq r_j$. \goodbreak 

\begin{theorem} 
Let $\phi$ be a symbol in $A (\D)$ whose image touches $\partial\D$ exactly at the points $\xi_1, \ldots, \xi_l$ and which is globally-regular. Then there are 
constants $\kappa$, $K$, $L > 0$, depending only on $\phi$, such that, for every $k \geq 1$: 
\begin{equation} \label{surprise} 
a_k (C_\phi) \leq K \, \bigg[ \frac{\omega^{- 1} (\kappa \, 2^{- N_k})}{\kappa \, 2^{- N_k}} \bigg]^{1/p} \,, 
\end{equation}
where $N_k$ is the largest integer such that $l N d_N < k$ and $d_N$ is the integer part of 
$\big[\log \frac{ \kappa \, 2^{- N}}{\omega^{- 1} (\kappa \, 2^{- N})} \big/ \log (\chi^{- p}) \big]+ 1$, with $0 < \chi < 1$ an absolute constant. 
\end{theorem} 

As a corollary, we get for lens maps $\lambda_\theta$ (as well as for polygonal maps), in the same way as Theorem~2.4 in \cite{estimates}, p.~550 (recall that then 
$\omega (h) \approx h^\theta$), the following upper bound. 
\begin{theorem} 
Let $\varphi = \lambda_\theta$ be the lens map of parameter $\theta$ acting on $H^p$, $1 < p < \infty$. Then, for positive constants $b$ and $c$ depending 
only on $\theta$ and $p$: 
\begin{displaymath} 
a_{n} (C_{\lambda_\theta}) \leq c \, \e^{- b \sqrt n} .
\end{displaymath} 
\end{theorem}

For the cusp map, we also have as in \cite{estimates}, Theorem~4.3 (here, $\omega(h) \approx 1 / \log (1/h)$). 

\begin{theorem} 
Let $\phi = \chi$ be the cusp map. For some positive constants $b$ and $c$ depending only on $p$, one has:
\begin{displaymath} 
a_n (C_\chi) \leq c\, \e^{- b \, n / \log n}  . 
\end{displaymath} 
\end{theorem} 

\goodbreak

\bigskip

\noindent
{\rm Daniel Li}, Univ Lille Nord de France, \\
U-Artois, Laboratoire de Math\'ematiques de Lens EA~2462 \\ 
\& F\'ed\'eration CNRS Nord-Pas-de-Calais FR~2956, \\
Facult\'e des Sciences Jean Perrin, Rue Jean Souvraz, S.P.\kern 1mm 18, \\
F-62\kern 1mm 300 LENS, FRANCE \\ 
daniel.li@euler.univ-artois.fr
\medskip

\noindent
{\rm Herv\'e Queff\'elec}, Univ Lille Nord de France, \\
USTL, Laboratoire Paul Painlev\'e U.M.R. CNRS 8524 \& 
F\'ed\'eration CNRS Nord-Pas-de-Calais FR~2956, \\
F-59\kern 1mm 655 VILLENEUVE D'ASCQ Cedex, 
FRANCE \\ 
Herve.Queffelec@univ-lille1.fr
\smallskip

\noindent
{\rm Luis Rodr{\'\i}guez-Piazza}, Universidad de Sevilla, \\
Facultad de Matem\'aticas, Departamento de An\'alisis Matem\'atico \& IMUS,\\ 
Apartado de Correos 1160,\\
41\kern 1mm 080 SEVILLA, SPAIN \\ 
piazza@us.es\par

\end{document}